\documentclass[12pt,a4paper]{amsart}
\usepackage{amsmath,amsfonts,amssymb,amsthm,xcolor}
\usepackage{tikz-cd}
\usepackage[pagebackref,colorlinks=true,linkcolor=blue,citecolor=red,urlcolor=blue]{hyperref}
\usepackage[a4paper,twoside,inner=99pt,outer=71pt,top=97pt,bottom=67pt,headheight=12pt,headsep=26pt,footskip=30pt]{geometry}
\newcommand{\R}{\mathbb R}
\newcommand{\Z}{\mathbb Z}
\newcommand{\N}{\mathbb N}
\newcommand{\F}{\mathbb F}
\newcommand{\cN}{\mathcal N}
\newcommand{\Heis}{\mathrm{Heis}}
\newcommand{\ProLieNil}{\mathsf{ProLieNil}}
\newcommand{\ProHeis}{\mathsf{ProHeis}}

\newtheorem{theorem}{Theorem}[section]
\newtheorem{corollary}[theorem]{Corollary}
\newtheorem{lemma}[theorem]{Lemma}

\theoremstyle{definition}
\newtheorem{definition}[theorem]{Definition}
\newtheorem{remark}[theorem]{Remark}
\newtheorem{example}[theorem]{Example}
\numberwithin{equation}{section}
\begin{document}
\title[Commuting Conditional Coproducts of Nilpotent Pro-Lie Groups]{Commuting Conditional Coproducts\\ of Nilpotent Pro-Lie Groups}
\author[F.G. Russo]{Francesco G. Russo}
\thanks{\\
Francesco G. Russo\newline
School of Science and Technology -- 
University of Camerino\newline
via Madonna delle Carceri 9, Camerino, Italy\newline
and\newline
Department of Mathematics and Applied Mathematics -- University of the Western Cape\newline
Private Bag X17, 7535, Bellville, South Africa\newline
and\newline
Department of Mathematics and Applied Mathematics -- University of Cape Town\newline
Private Bag X1, Rondebosch 7701, Cape Town, South Africa}
\begin{abstract}
In recent works of the literature, the special behaviour of coproducts in the categories of profinite abelian groups, locally compact abelian groups and abelian pro-Lie groups has been investigated. In the category of abelian pro-Lie groups one can formulate the notion of conditional coproduct, which agrees with that of cartesian product. An obstruction arises when we consider nonabelian topological groups, since the images of distinct factors of a cartesian product must commute. The notion of commuting conditional coproduct allows us to extend the abelian result to nonabelian pro-Lie groups, with particular attention to nilpotent pro-Lie groups of class at most two. We also describe the behaviour of semidirect products under cartesian products and recall a Fitting-type lemma for pro-Lie groups.
\medskip

\noindent\textit{Keywords:} Pro-Lie groups, nilpotent groups, coproducts, Heisenberg groups.

\noindent\textit{Mathematics Subject Classification 2020:} 22E25, 22D05, 22A05; 20J15, 43A95.
\end{abstract}
\maketitle

\section{A brief introduction}

We will consider only Hausdorff topological groups. Assuming some familiarity with the terminology and notations of \cite{hn,hofmor,prolie,hannaneumann}, we generalize \cite[Theorem 4]{hhrneatnote} to nonabelian pro-Lie groups, in the context of the notion of varieties of groups of Morris and others \cite{gm,hofmorstr,morris3,morris4,morris5}. We introduce new concepts in Definitions~\ref{d:convergent} and \ref{d:condprod}, looking at certain families of Heisenberg groups as a source of inspiration.

Theorem~\ref{t:tychonovnonabelian} is our first main result and generalizes \cite[Theorem 4]{hhrneatnote} to pro-Lie groups of nilpotence class at most two. Its proof in fact applies to pro-Lie groups of arbitrary nilpotence class, and even without a nilpotence assumption. Theorem~\ref{t:products-semidirect} describes the corresponding behaviour of semidirect products for arbitrary topological groups. The constructions of Example~\ref{Heisenberg} and Corollary~\ref{c:proheis-conditional} illustrate these results for Heisenberg groups, which are important in pure and applied mathematics. Finally a Fitting-type lemma is given by Theorem~\ref{t:fitting}, in order to recall how the nilpotence class of a closed subgroup generated by two normal nilpotent subgroups may be controlled. The commuting condition is essential for the universal property of a cartesian product in a nonabelian category.

After the present introduction in Section~1, we have Section~2, which illustrates some concepts taken from the existing literature and provides the background for the new notions. In Section~3 we give details on the example which suggested the main results of the present paper, namely Example~\ref{Heisenberg}. Then we prove our results in Section~4, indeed Theorems~\ref{t:tychonovnonabelian}, \ref{t:products-semidirect} and \ref{t:fitting}.

\section{Varieties of nilpotent pro-Lie groups}

Regarding the notion of compact Lie group, there are different ways to introduce it, but we may look at one of the equivalent conditions in \cite[Corollary 2.40, Definition 2.41]{hofmor}. In particular, we may think of a compact Lie group as a compact group possessing a faithful continuous finite dimensional orthogonal (or unitary) representation. Concerning noncompact Lie groups, one should introduce the notion of finite dimensional real Lie group and so one can proceed as in \cite[Definition 2.1, Proposition 2.2]{prolie}. See \cite{hn,hofmor,prolie} for basic theory of Lie groups and pro-Lie groups.

\begin{definition}[See \cite{prolie}, Definition 3.25]\label{preciseprolie}
A topological group $G$ has arbitrarily small co-Lie subgroups if every identity neighbourhood of $G$ contains a closed normal subgroup $N$ such that $G/N$ is a Lie group. We say that $G$ is a pro-Lie group if $G$ has arbitrarily small co-Lie subgroups and is complete with respect to the two-sided uniformity.
\end{definition}

For a pro-Lie group $G$, the family
\begin{equation}\label{eq:colie}
\cN(G)=\{N\trianglelefteq G\mid N\text{ is closed and }G/N\text{ is a Lie group}\}
\end{equation}
is a filter basis of closed normal subgroups converging to the trivial subgroup, and $G\simeq\varprojlim_{N\in\cN(G)}G/N$ as a topological group; see \cite[Theorem 3.27]{prolie}, or \cite[pp. 63--68]{hofmorstr}. A topological group $G$ is almost connected if the factor group $G/G_0$ of $G$ modulo the connected component $G_0$ of the identity is compact. See also \cite[Definitions A, B, C, Chapter 3]{prolie}. Among almost connected groups we may find compact groups and connected locally compact groups. Every almost connected locally compact group is a pro-Lie group. Locally compact abelian groups are abelian pro-Lie groups, see \cite[Example 14.54]{prolie}.

\begin{definition}[See \cite{prolie}, pp. 644--646]\label{tricky}
A topological group $G$ is called
\begin{itemize}
\item[(i).] \textit{Polish} if its topology is separable and completely metrizable;
\item[(ii).] $\sigma$-\textit{compact} if it is a countable union of compact subspaces;
\item[(iii).] \textit{Compactly generated} if there is a compact subset $K$ of $G$ such that $\langle K\rangle=G$;
\item[(iv).] \textit{Compactly topologically generated} if there is a compact subset $K$ such that $\overline{\langle K\rangle}=G$.
\end{itemize}
\end{definition}
It is not immediate whether groups satisfying the conditions in Definition~\ref{tricky} are pro-Lie groups or not.

\begin{remark}\label{sigmacompactness}
According to \cite[Remark 5.22 and Appendix 2]{prolie}, we have that
\begin{itemize}
\item[(i).] Every almost connected locally compact group is compactly generated.
\item[(ii).] Every compactly generated group is $\sigma$-compact.
\item[(iii).] If $d$ and $m$ are two nonnegative integers, $K$ a compact abelian group, $\R^d$ the direct sum of $d$ copies of the additive group $\R$ of the reals, $\Z^m$ the direct sum of $m$ copies of the additive group $\Z$ of the integers, then \cite[Theorem 7.57]{hofmor} shows that every compactly generated locally compact abelian group is isomorphic to a direct sum $\R^d\oplus\Z^m\oplus K$, which is an example of an abelian pro-Lie group.
\item[(iv).] A topological group whose underlying topological space is $\sigma$-compact and a Baire space is a locally compact group.
\item[(v).] A $\sigma$-compact Polish group is locally compact.
\item[(vi).] A compactly generated topological group whose underlying topological space is a Baire space is locally compact.
\end{itemize}
\end{remark}

In \cite{hn,hofmor,prolie} there are several examples which describe the conditions in Definition~\ref{tricky} and Remark~\ref{sigmacompactness}. For instance, finite abelian groups satisfy trivially all the conditions in Definition~\ref{tricky}; see \cite[Chapters 5, 9, 14]{prolie}. Already cartesian products of abelian groups show that one can keep the condition of being a pro-Lie group and lose that of being a locally compact group.

\begin{example}\label{topologyofprolie}
Referring to \cite[Remark 5.22]{prolie}, the cartesian product
\begin{equation}\label{eq:Rprod}
\R^\N=\{(g_n)_{n\in\N}\mid g_n\in\R\}
\end{equation}
of countably many copies of $\R$ with its usual topology is a topological group if we define pointwise addition and consider the product topology on it. This is an example of a Polish group, which is neither discrete nor locally compact. Note also that $\R^\N$ is a pro-Lie group. One can also check that its closed subgroup
\begin{equation}\label{eq:Zprod}
\Z^\N=\{(g_n)_{n\in\N}\mid g_n\in\Z\}
\end{equation}
is a Polish group which is not locally compact. Therefore already in the abelian case we have examples for Definitions~\ref{preciseprolie} and \ref{tricky}.
\end{example}

For commutators we use $[x,y]=x^{-1}y^{-1}xy$. The general term of the (abstract) lower central series of a group $G$ will be denoted by $\gamma_i(G)$ and the general term of the upper central series of $G$ by $Z_i(G)$, where $i\in\N$. While $Z_i(G)$ is always a closed subgroup of a Hausdorff topological group $G$, this is no longer true for $\gamma_i(G)$; the corresponding closed lower central series has terms $\overline{\gamma_i(G)}$. See details in \cite{hofmor} and \cite[Chapter 10]{prolie}. We say that $G$ is nilpotent of class at most $c$ if $\gamma_{c+1}(G)=1$, or equivalently $Z_c(G)=G$. The least such $c$ is the nilpotence class $c(G)$ of $G$. These definitions are unchanged if one takes the closure of each term of the lower central series. Nilpotent groups help to find large classes of topological groups which generalize abelian pro-Lie groups. Following Morris \cite{morris3,morris4,morris5} and Neumann \cite{hannaneumann}, it is appropriate to introduce $\overline{\mathbf S}$, $\overline{\mathbf Q}$, $\mathbf C$ on a class of topological groups to be, respectively, the operations of taking arbitrary closed subgroups, quotients by closed normal subgroups and arbitrary cartesian products.

\begin{definition}[See \cite{morris3,morris4,morris5}]\label{Hausdorffvariety}
A nonempty class of topological groups $\Omega$, closed under topological isomorphisms, is called a variety of groups if the following conditions are satisfied:
\begin{itemize}
\item[(a).] If $G\in\Omega$ and $H$ is a closed subgroup of $G$, then $H\in\Omega$;
\item[(b).] If $N$ is a closed normal subgroup of $G$ and $G\in\Omega$, then $G/N\in\Omega$;
\item[(c).] If $(G_i)_{i\in I}$ is a family of groups in $\Omega$, then $\prod_{i\in I}G_i\in\Omega$.
\end{itemize}
\end{definition}

We may rewrite (a) above as $\overline{\mathbf S}\Omega=\Omega$, (b) as $\overline{\mathbf Q}\Omega=\Omega$ and (c) as $\mathbf C\Omega=\Omega$. With the help of $\overline{\mathbf S}$, $\overline{\mathbf Q}$ and $\mathbf C$, it is possible to construct the variety generated by a given class of topological groups by taking the smallest class containing it which is closed under the three operations. See \cite[Theorem 3.3, Examples 3.4, and Chapter 4]{prolie} for the relevance of these operations to projective limits of Lie groups, and for the particular difficulties concerning quotients of pro-Lie groups.

\begin{lemma}\label{nilpotenclasstwovariety}
For every nonnegative integer $c$, the class
\[
\ProLieNil_c=\{G\mid G\text{ is a pro-Lie group and }\gamma_{c+1}(G)=1\}
\]
is closed under taking closed subgroups, arbitrary cartesian products and quotients by closed normal subgroups whenever the quotient is complete.
\end{lemma}

\begin{proof}
Let $G\in\ProLieNil_c$. First, let $H$ be a closed subgroup of $G$. By \cite[Proposition 3.2]{prolie}, every closed subgroup of a pro-Lie group is again a pro-Lie group. Moreover, $\gamma_i(H)\leq\gamma_i(G)$ for every $i\geq1$. Hence $\gamma_{c+1}(H)\leq\gamma_{c+1}(G)=1$. Therefore $H$ is in $\ProLieNil_c$.

Secondly, let $N$ be a closed normal subgroup of $G$ and suppose that $G/N$ is complete. By \cite[Theorem 4.1(i)]{prolie}, the quotient $G/N$ is again a pro-Lie group. If $\pi\colon G\to G/N$ denotes the quotient homomorphism, then $\gamma_i(G/N)=\pi(\gamma_i(G))=\gamma_i(G)N/N$ for every $i\geq1$. Consequently, $\gamma_{c+1}(G/N)=\gamma_{c+1}(G)N/N=1$ and so $G/N$ is in $\ProLieNil_c$.

Finally, let $(G_j)_{j\in J}$ be a family of groups in $\ProLieNil_c$, and consider the cartesian product
\begin{equation}\label{eq:genericprod}
\prod_{j\in J}G_j=\{(g_j)_{j\in J}\mid g_j\in G_j\}
\end{equation}
with the product topology. Arbitrary cartesian products of pro-Lie groups are again pro-Lie groups. Furthermore, commutators in \eqref{eq:genericprod} are computed pointwise. Therefore, for every $i\geq1$,
\begin{equation}\label{eq:gammaproduct}
\gamma_i\left(\prod_{j\in J}G_j\right)\subseteq\prod_{j\in J}\gamma_i(G_j).
\end{equation}
In particular,
\begin{equation}\label{eq:nilprod}
\gamma_{c+1}\left(\prod_{j\in J}G_j\right)\subseteq\prod_{j\in J}\gamma_{c+1}(G_j)=\{1\}.
\end{equation}
This allows us to conclude that \eqref{eq:genericprod} belongs to $\ProLieNil_c$, and the result follows.
\end{proof}

The completeness hypothesis in the quotient assertion cannot be omitted, even for abelian pro-Lie groups: \cite[Corollary 4.11]{prolie} gives a closed subgroup of $\R^{[0,1]}$ with incomplete quotient. Thus $\ProLieNil_c$ need not be a variety in the sense of Definition~\ref{Hausdorffvariety}; indeed it is not such a variety for any $c\geq1$. In the following we regard $\ProLieNil_c$ as a full subcategory of the category of Hausdorff topological groups, with continuous homomorphisms as morphisms. Closure under arbitrary products and closed subgroups is all that will be needed.

\begin{definition}\label{heisenbergvariety}
A variety of groups in Definition~\ref{Hausdorffvariety} is called a variety of nilpotent pro-Lie groups if all the members of $\Omega$ are nilpotent pro-Lie groups.
\end{definition}

The consideration of compact abelian groups shows that there exist nontrivial varieties of nilpotent pro-Lie groups. In fact, coproducts can be very interesting to describe for pro-Lie groups and profinite groups, see \cite{chir1,chir2,gm,hhrneat,hhrneatnote,hm77,Neukirch71}.

\section{Some considerations on pro-Heisenberg groups}

We may begin by considering the following example, whose finite dimensional description is recalled from \cite[Example 4.4]{br}, in order to extend the ideas in \cite{hhrneatnote} to nilpotent pro-Lie groups which are not necessarily abelian. These groups are largely studied in mathematical physics, see also \cite{rw}, and in a different perspective in \cite[Example 14.23]{prolie}, where pro-Heisenberg groups over weakly complete topological spaces are considered.

\begin{example}\label{Heisenberg}
Consider the three dimensional connected Heisenberg group on $\R$
\begin{equation}\label{largeheisenberg1}
\Heis(\R)=\left\{\begin{pmatrix}1&a&c\\0&1&b\\0&0&1\end{pmatrix}\ \middle|\ a,b,c\in\R\right\}
\end{equation}
represented as a group of matrices in the general linear group $\mathrm{GL}(3,\R)$. This is a locally compact nilpotent pro-Lie group, see \cite[Examples 10.44 and 14.39]{prolie}. Its structural properties are well known, see \cite{br,hn,hofmor,prolie,rw}. In particular, $\Heis(\R)$ has one dimensional center and two dimensional central quotient, i.e.:
\begin{equation}\label{largeheisenberg3}
Z(\Heis(\R))=[\Heis(\R),\Heis(\R)]\simeq\R
\quad\text{and}\quad \Heis(\R)/Z(\Heis(\R))\simeq\R^2.
\end{equation}
Moreover it can be written as a topological semidirect product in the following form
\begin{equation}\label{splitting}
\Heis(\R)=H\rtimes K\simeq\R^2\rtimes\R,
\end{equation}
where
\begin{equation}\label{AB}
\begin{aligned}
H&=\left\{\begin{pmatrix}1&a&c\\0&1&0\\0&0&1\end{pmatrix}\ \middle|\ a,c\in\R\right\},\\
K&=\left\{\begin{pmatrix}1&0&0\\0&1&b\\0&0&1\end{pmatrix}\ \middle|\ b\in\R\right\}.
\end{aligned}
\end{equation}
Here the action of $K\simeq\R$ on $H\simeq\R^2$ is given, in the coordinates of \eqref{AB}, by $b\cdot(a,c)=(a,c-ab)$. Note that $Z_2(\Heis(\R))=\Heis(\R)$, that is, $\Heis(\R)$ has nilpotence class two. Moreover \eqref{largeheisenberg1} possesses a normal vector subgroup, namely its center, but does not split over it: such a splitting would be a direct product of the center with the abelian central quotient, contradicting noncommutativity. See \cite[Lemma 11.13, Corollary 11.14]{prolie}.

Now it is possible to form, with
\begin{equation}\label{uniformization}
A_1=A_2=A_3=\cdots=A_n=A_{n+1}=\cdots=\Heis(\R),
\end{equation}
the cartesian product
\begin{equation}\label{proheis}
P=\prod_{n\in\N}A_n=\Heis(\R)^\N=\{(g_n)_{n\in\N}\mid g_n\in\Heis(\R)\},
\end{equation}
which has
\begin{equation}\label{eq:heiscenter}
Z(P)=Z\left(\prod_{n\in\N}A_n\right)=\prod_{n\in\N}Z(A_n)=Z(\Heis(\R))^\N\simeq\R^\N,
\end{equation}
where the equality follows from the coordinatewise group operation. In this particular case, also
\begin{equation}\label{eq:heisderived}
[P,P]=[\Heis(\R),\Heis(\R)]^\N=Z(\Heis(\R))^\N\simeq\R^\N,
\end{equation}
since every element of $Z(\Heis(\R))$ is a single commutator. Indeed, the commutator of the matrices with entries $a=1$, $b=c=0$ and $b=t$, $a=c=0$, respectively, is the central matrix with entry $c=t$. Choosing these matrices coordinatewise proves \eqref{eq:heisderived}. Note that
\begin{equation}\label{eq:heisab}
P/[P,P]\simeq(\R\times\R)^\N=\{(a_n,b_n)_{n\in\N}\mid a_n,b_n\in\R\}.
\end{equation}
We may also consider the (restricted) direct sum
\begin{equation}\label{proheissum}
\begin{aligned}
Q&=\bigoplus_{n\in\N}A_n=\Heis(\R)^{(\N)}\\
 &=\{(g_n)_{n\in\N}\mid g_n\in\Heis(\R),\ g_n=1\text{ for all but finitely many }n\},
\end{aligned}
\end{equation}
with the subspace topology inherited from $P$, and get in analogy
\begin{equation}\label{eq:heissumcenter}
Z(Q)=Z(\Heis(\R))^{(\N)}\simeq\R^{(\N)},
\end{equation}
\begin{equation}\label{eq:heissumderived}
[Q,Q]=[\Heis(\R),\Heis(\R)]^{(\N)}=Z(\Heis(\R))^{(\N)}\simeq\R^{(\N)}.
\end{equation}
The same commutator calculation proves \eqref{eq:heissumderived}, using matrices with finite support. From Lemma~\ref{nilpotenclasstwovariety}, $P$ belongs to $\ProLieNil_2$. It is a Polish nilpotent pro-Lie group which is not locally compact. On the other hand, $Q$ is a proper dense subgroup of $P$ and hence is not complete with the subspace topology. In particular, $Q$ is not a pro-Lie group with this topology. The direct sum in \eqref{proheissum} is not to be identified with the categorical coproduct in $\ProLieNil_2$.

The subgroup $H$ in \eqref{AB} of $\Heis(\R)$ gives rise to a corresponding subgroup $H_n\simeq\R^2$ in $A_n$ as in \eqref{uniformization}, and similarly for the copies $K_n$ of $K$. Therefore we get
\begin{equation}\label{eq:Hprod}
\prod_{n\in\N}H_n=H^\N=\{(z_n,u_n)_{n\in\N}\mid(z_n,u_n)\in H\}\simeq(\R\times\R)^\N,
\end{equation}
\begin{equation}\label{eq:Kprod}
\prod_{n\in\N}K_n=K^\N=\{(v_n)_{n\in\N}\mid v_n\in K\}\simeq\R^\N
\end{equation}
and we find a decomposition in topological semidirect product
\begin{equation}\label{sp}
\Heis(\R)^\N=H^\N\rtimes K^\N.
\end{equation}
On the other hand, we can work at the level of finite-support subgroups, finding that
\begin{equation}\label{eq:HKsum}
\bigoplus_{n\in\N}H_n=H^{(\N)}\simeq(\R\times\R)^{(\N)}
\quad\text{and}\quad\bigoplus_{n\in\N}K_n=K^{(\N)}\simeq\R^{(\N)}
\end{equation}
and so
\begin{equation}\label{spcop}
\Heis(\R)^{(\N)}=H^{(\N)}\rtimes K^{(\N)},
\end{equation}
where all finite-support groups carry their respective subspace topologies. The actions in \eqref{sp} and \eqref{spcop} are coordinatewise.

There are several consequences for \eqref{sp} and \eqref{spcop}. Firstly, \eqref{sp} shows that for a given $(g_n)_{n\in\N}\in\Heis(\R)^\N$ there exist unique $(h_n)_{n\in\N}\in H^\N$ and $(k_n)_{n\in\N}\in K^\N$ such that $g_n=h_nk_n$ for every $n$. This corresponds also at the level of \eqref{spcop} with the additional condition that $h_n=k_n=1$ for all but finitely many $n$.

Secondly, write $H_{[m]}=\prod_{n=1}^m H_n$, $K_{[m]}=\prod_{n=1}^m K_n$ and $P_{[m]}=\prod_{n=1}^m A_n=H_{[m]}\rtimes K_{[m]}$. With the coordinate projections as connecting morphisms, which preserve the semidirect decompositions and their actions, we have
\begin{equation}\label{eq:heislim}
\varprojlim_m P_{[m]}\simeq\left(\varprojlim_m H_{[m]}\right)\rtimes
\left(\varprojlim_m K_{[m]}\right)\simeq P.
\end{equation}
Using instead the embeddings obtained by adjoining identity coordinates, the corresponding algebraic direct limits give $Q$ and its factors in \eqref{spcop}. We equip these algebraic direct limits with the subspace topologies specified above, without asserting a colimit property in the category of Hausdorff topological groups. While the factors $H$ and $K$ are abelian, Pontryagin duality of abelian groups does not identify the nonabelian direct sum $Q$ with a categorical coproduct. The conditional universal property of $P$ is the subject of Section~4.
\end{example}

It is useful to mention here an interesting situation for finite groups, in order to show the relevance of a splitting in Example~\ref{Heisenberg}.

\begin{example}\label{hamiltonian}
The quaternion group $Q_8$ has order eight and is a $2$-group of nilpotence class two. There are no nontrivial subgroups $H$ and $K$ in $Q_8$ such that $Q_8=H\rtimes K$. In fact, every nontrivial subgroup of $Q_8$ contains its unique involution, so two such subgroups cannot have trivial intersection. Thus a decomposition as in \eqref{splitting} need not exist even for compact groups of nilpotence class two. In such cases one cannot use that decomposition to describe the product, although the universal property in Theorem~\ref{t:tychonovnonabelian} below still applies.
\end{example}

\section{Splitting of nilpotent pro-Lie groups}

We want to design an appropriate context in order to explore what we have just noted for the pro-Heisenberg group in Example~\ref{Heisenberg}, generalizing \cite[Definition 1]{hhrneatnote}.

\begin{definition}\label{d:convergent}
Consider $\ProLieNil_2$ as in Lemma~\ref{nilpotenclasstwovariety}. Let $(A_j)_{j\in J}$ be a family of objects in $\ProLieNil_2$ with $J$ an arbitrary set of indices, $H$ also an object in $\ProLieNil_2$, and $\F=(\phi_j\colon A_j\to H)_{j\in J}$ a family of continuous homomorphisms. We say that $\F$
\begin{itemize}
\item[(a)] is \emph{convergent} if for every neighbourhood $U$ of the identity of $H$
\[
J_U:=\{j\in J\mid\phi_j(A_j)\not\subseteq U\}\text{ is finite};
\]
\item[(b)] is \emph{commuting and convergent} if it is convergent and
\[
[\phi_i(A_i),\phi_j(A_j)]=1\qquad\text{whenever }i\ne j.
\]
\end{itemize}
\end{definition}

Of course, $\ProLieNil_1$ is exactly the category of abelian pro-Lie groups with the corresponding morphisms and in this case there is no difference between the conditions (a) and (b) in Definition~\ref{d:convergent}, since every family of morphisms in $\ProLieNil_1$ which is convergent is also commuting, but this is no longer true in $\ProLieNil_2$ and more generally in $\ProLieNil_c$ when $c\geq2$. The same definition makes sense for arbitrary Hausdorff topological groups.

\begin{example}\label{ex:P-convergent}
We may adapt \cite[Example 2]{hhrneatnote} to the present context. For any family $(A_j)_{j\in J}$ of pro-Lie groups, not necessarily abelian, let $H=\prod_{j\in J}A_j$ be the cartesian product with the product topology. Then the family $\F=(\tau_j)_{j\in J}$ of canonical embeddings $\tau_j\colon A_j\to H$, given by $\tau_j(a)=(a_l)_{l\in J}$ where $a_l=a$ if $l=j$ and $a_l=1$ otherwise, is commuting and convergent. In fact, each identity neighbourhood in $H$ contains a basic product neighbourhood restricting only finitely many coordinates, so it contains the images of all but finitely many factors. In particular, this applies to Example~\ref{Heisenberg} with $J=\N$, $A_j=\Heis(\R)$ and $H=\Heis(\R)^\N$.
\end{example}

We define the commuting conditional coproduct by means of a universal property, resembling the  conditional coproduct of abelian groups in \cite[Definition 3]{hhrneatnote}.

\begin{definition}\label{d:condprod}
Let $(A_j)_{j\in J}$ be a family of objects in $\ProLieNil_2$. A pair $(G,(\tau_j)_{j\in J})$, where $G\in\ProLieNil_2$ and the family of continuous homomorphisms $\tau_j\colon A_j\longrightarrow G$ is commuting and convergent, is called a commuting conditional coproduct of $(A_j)_{j\in J}$ if, for every object $H$ in $\ProLieNil_2$ and every commuting and convergent family of continuous homomorphisms $\psi_j\colon A_j\longrightarrow H$, there exists a unique continuous homomorphism $\omega\colon G\longrightarrow H$ such that $\psi_j=\omega\circ\tau_j$ for every $j\in J$.
\end{definition}

In the case of $\ProLieNil_1$ (with the definition taken in that category), Definition~\ref{d:condprod} agrees with \cite[Definition 3]{hhrneatnote}, because the two conditions of Definition~\ref{d:convergent} agree. In particular, for any $j\in J$, we can visualize the universal property of Definition~\ref{d:condprod} with the following commutative diagram:
\begin{equation}\label{ccpd}
\begin{tikzcd}[column sep=large,row sep=large]
A_j \arrow[r,"\tau_j"] \arrow[dr,"\psi_j"'] & G \arrow[d,dashed,"\exists!\,\omega"]\\
 & H
\end{tikzcd}
\end{equation}
The same $\omega$ must work simultaneously for all $j\in J$. The universal property determines a commuting conditional coproduct uniquely up to a unique isomorphism preserving the canonical morphisms: the universal properties of two such pairs give morphisms in both directions, and their composites must be the respective identities by uniqueness.

\begin{definition}\label{proheiscategory}
Consider the full subcategory $\ProHeis$ of $\ProLieNil_2$, whose objects are the groups admitting a topological semidirect product decomposition
\[
M=E\rtimes F,
\]
where $E$ is a closed normal abelian pro-Lie subgroup of $M$ and $F$ is a closed abelian pro-Lie subgroup of $M$. In other words, the multiplication map $E\times F\to M$ is a homeomorphism and $E\cap F=\{1\}$.
\end{definition}

The groups $\Heis(\R)^n$ ($n\in\N$) and $\Heis(\R)^\N$ belong to $\ProHeis$. Example~\ref{hamiltonian} shows that there are groups belonging to $\ProLieNil_2$ but neither to $\ProHeis$ nor to $\ProLieNil_1$. Abelian pro-Lie groups belong to $\ProHeis$ by taking the second factor to be trivial.

\begin{theorem}[First Main Theorem]\label{t:tychonovnonabelian}
If $(A_j)_{j\in J}$ is a family of groups in $\ProLieNil_2$, then the cartesian product
\[
P=\prod_{j\in J}A_j,
\]
together with the canonical embeddings $\tau_j\colon A_j\longrightarrow P$, is a commuting conditional coproduct of $(A_j)_{j\in J}$ in $\ProLieNil_2$.
\end{theorem}

\begin{proof}
Since the category of pro-Lie groups is closed under the formation of arbitrary products by \cite[Theorem 3.3]{prolie}, $P$ is a pro-Lie group, and commutators in $P$ are computed pointwise. Moreover, as each $A_j$ has nilpotence class at most two, so does $P$. Thus $P$ belongs to $\ProLieNil_2$. It remains to check the conditions in Definition~\ref{d:condprod}. For distinct indices $i,j\in J$, the subgroups $\tau_i(A_i)$ and $\tau_j(A_j)$ commute, because they have disjoint supports. The family $(\tau_j)_{j\in J}$ is also convergent. Indeed, if $U$ is an identity neighbourhood in $P$, then $U$ contains a basic identity neighbourhood
\begin{equation}\label{eq:basic}
V=\prod_{j\in J}V_j,
\end{equation}
where $V_j=A_j$ for all but finitely many $j$. Hence
\begin{equation}\label{eq:tail}
\tau_j(A_j)\subseteq V\subseteq U
\end{equation}
for all but finitely many $j$.

Now let $H$ be another group belonging to $\ProLieNil_2$, and let $\psi_j\colon A_j\longrightarrow H$ ($j\in J$) be a commuting and convergent family of continuous homomorphisms. We follow the method of \cite[Lemma 6 and proof of Theorem 4]{hhrneatnote}, keeping track of the commutativity of distinct images. Consider the filter basis $\cN(H)$ of all closed normal subgroups $N$ of $H$ for which $H/N$ is a finite dimensional Lie group. Since $H$ is a pro-Lie group, we have $H\simeq\varprojlim_{N\in\cN(H)}H/N$ by \cite[Theorem 3.27]{prolie}. Fix $N\in\cN(H)$. Since the Lie group $H/N$ has no small subgroups, it has an identity neighbourhood $V_N$ containing no nontrivial subgroup. If $q_N\colon H\longrightarrow H/N$ is the quotient map, then $U_N=q_N^{-1}(V_N)$ is an identity neighbourhood in $H$. By convergence, only finitely many $j$ have $\psi_j(A_j)\not\subseteq U_N$. If $\psi_j(A_j)\subseteq U_N$, then $q_N(\psi_j(A_j))$ is a subgroup contained in $V_N$, and therefore $q_N(\psi_j(A_j))=1$. Consequently, as in \cite[Lemma 6]{hhrneatnote}, the set
\begin{equation}\label{eq:JN}
J_N=\{j\in J\mid\psi_j(A_j)\not\subseteq N\}
\end{equation}
is finite, and
\begin{equation}\label{eq:tailN}
\psi_j(A_j)\subseteq N\qquad\text{for every }j\notin J_N.
\end{equation}

This is the step which helps to construct $\omega$. In fact we may define first
\begin{equation}\label{omegaN}
\omega_N\colon P\longrightarrow H/N,\qquad x=(a_j)_{j\in J}\longmapsto\prod_{j\in J_N}q_N(\psi_j(a_j)),
\end{equation}
which gives a finite product and is independent of the order chosen on $J_N$, because the subgroups $\psi_i(A_i)$ and $\psi_j(A_j)$ commute whenever $i\ne j$. Then we observe that the map $\omega_N$ is a continuous homomorphism. Indeed, it depends on only finitely many coordinates, and for $x=(a_j)$ and $y=(b_j)$ one has
\begin{equation}\label{eq:omegahom}
\begin{aligned}
\omega_N(xy)&=\prod_{j\in J_N}q_N(\psi_j(a_jb_j))\\
 &=\prod_{j\in J_N}\bigl(q_N(\psi_j(a_j))q_N(\psi_j(b_j))\bigr)\\
 &=\left(\prod_{j\in J_N}q_N(\psi_j(a_j))\right)
 \left(\prod_{j\in J_N}q_N(\psi_j(b_j))\right)=\omega_N(x)\omega_N(y),
\end{aligned}
\end{equation}
where the third equality uses the pairwise commutativity of the images belonging to distinct indices. At this point, we may construct the requested $\omega$ via the $\omega_N$ above and using the structure of projective limit. If $M,N\in\cN(H)$ and $M\subseteq N$, then $J_N\subseteq J_M$, and the canonical quotient homomorphism
\begin{equation}\label{eq:piMN}
\pi_{MN}\colon H/M\longrightarrow H/N
\end{equation}
satisfies
\begin{equation}\label{omegaMN}
\omega_N=\pi_{MN}\circ\omega_M.
\end{equation}
Indeed, all terms corresponding to indices in $J_M\setminus J_N$ vanish modulo $N$. Thus $(\omega_N)_{N\in\cN(H)}$ is a compatible family of continuous homomorphisms. By the universal property of the inverse limit \cite[Lemma 1.26]{hofmor}, there exists a unique continuous homomorphism
\begin{equation}\label{eq:omega}
\omega\colon P\longrightarrow H
\end{equation}
such that
\begin{equation}\label{eq:qomega}
q_N\circ\omega=\omega_N\qquad\text{for every }N\in\cN(H).
\end{equation}
This is the required $\omega$; in fact, for $a\in A_j$, the definitions imply $q_N(\omega(\tau_j(a)))=q_N(\psi_j(a))$ for every $N\in\cN(H)$, and $\bigcap_{N\in\cN(H)}N=\{1\}$ gives $\omega\circ\tau_j=\psi_j$ for all $j\in J$. It remains to prove uniqueness. Let
\begin{equation}\label{eq:Pfin}
P_{\rm fin}=\{(a_j)_{j\in J}\in P\mid a_j=1\text{ for all but finitely many }j\}.
\end{equation}
This subgroup is dense in $P$ and every element of $P_{\rm fin}$ is a finite product of elements belonging to the canonical subgroups $\tau_j(A_j)$. Therefore any continuous homomorphism $\eta\colon P\to H$ satisfying $\eta\circ\tau_j=\psi_j$ for every $j$ must agree with $\omega$ on $P_{\rm fin}$. Since $\eta$ and $\omega$ are continuous and $H$ is Hausdorff, they agree on all of $P$. Hence $(P,(\tau_j)_{j\in J})$ has the asserted universal property, that is, Definition~\ref{d:condprod} is satisfied, and the proof follows.
\end{proof}

The proof of Theorem~\ref{t:tychonovnonabelian} does not use nilpotence, apart from establishing that $P$ belongs to the prescribed category. In fact, if $(A_j)_{j\in J}$ is a family of Hausdorff topological groups and $H$ is a pro-Lie group, a family of continuous homomorphisms $\psi_j\colon A_j\to H$ extends uniquely to a continuous homomorphism $\prod_{j\in J}A_j\to H$ if and only if it is commuting and convergent. Sufficiency follows from the same proof starting with \eqref{eq:JN}, and necessity follows by composition with the canonical embeddings, which are commuting and convergent. Thus the same conditional coproduct theorem holds in the category of all pro-Lie groups and in any full subcategory of it closed under arbitrary products, in particular in $\ProLieNil_c$ for every finite $c$.

\begin{corollary}\label{c:proheis-conditional}
Let $(A_j)_{j\in J}$ be a family of groups in $\ProHeis$, such that $A_j=E_j\rtimes F_j$ are topological semidirect products. Then their commuting conditional coproduct, both in $\ProHeis$ and in $\ProLieNil_2$, is
\[
\prod_{j\in J}A_j\simeq\left(\prod_{j\in J}E_j\right)\rtimes\left(\prod_{j\in J}F_j\right),
\]
with the coordinatewise action.
\end{corollary}
\begin{proof}
Regrouping the coordinates of $\prod_{j\in J}(E_j\rtimes F_j)$ gives the indicated topological semidirect product, by the same coordinatewise argument as in Theorem~\ref{t:products-semidirect} below. Its factors are abelian pro-Lie groups, so the product belongs to $\ProHeis$. The universal property follows from Theorem~\ref{t:tychonovnonabelian} and the fact that $\ProHeis$ is a full subcategory of $\ProLieNil_2$.
\end{proof}

We may describe this behaviour in a more general setting.

\begin{theorem}[Second Main Theorem]\label{t:products-semidirect}
Let $G=H\rtimes K$ be a topological semidirect product, where $H$ is a closed normal subgroup of $G$ and $K$ is a closed subgroup (not necessarily normal) acting on $H$ by automorphisms via a jointly continuous action. Then for every index set $J$
\[
G^J\simeq H^J\rtimes K^J,
\]
where the action of $K^J$ on $H^J$ is coordinatewise. Moreover, if $G^{(J)}$, $H^{(J)}$ and $K^{(J)}$ denote the respective finite-support subgroups equipped with the subspace topologies of the cartesian powers, then
\[
G^{(J)}\simeq H^{(J)}\rtimes K^{(J)}.
\]
\end{theorem}

\begin{proof}
Assume $G=H\rtimes_\alpha K$ is given by the action $\alpha$ as in the assumptions. Then $G^J$ is naturally identified with $(H\rtimes_\alpha K)^J$ and every element can be written uniquely as $(g_j)_{j\in J}=(h_jk_j)_{j\in J}$ with $h_j\in H$ and $k_j\in K$. Define
\begin{equation}\label{eq:Phi}
\begin{aligned}
\Phi\colon G^J&\longrightarrow H^J\rtimes_{\alpha^J} K^J,\\
(h_jk_j)_{j\in J}&\longmapsto\bigl((h_j)_{j\in J},(k_j)_{j\in J}\bigr).
\end{aligned}
\end{equation}
The action is coordinatewise, $\alpha^J((k_j))((h_j))=(\alpha(k_j)(h_j))$, and the multiplication is given by $(h_j,k_j)(h'_j,k'_j)=(h_j\alpha(k_j)(h'_j),k_jk'_j)$ as in the definition of the external semidirect product. Thus $\Phi$ is a homomorphism. The product of the homeomorphisms defining the topological semidirect product, followed by regrouping of the coordinates, shows that $\Phi$ and its inverse are continuous. Hence $\Phi$ is a topological group isomorphism and the first part of the theorem follows.

The same map restricts to the finite-support subgroups, since $h_jk_j=1$ if and only if $h_j=k_j=1$. Its restriction and the inverse restriction are continuous in the subspace topologies, which proves the second part.
\end{proof}

We shall mention that a classical result of Hans Fitting (see \cite{huppert} and \cite[Theorem 5.2.8]{robinson}) shows that the product of two normal nilpotent subgroups $H$ and $K$ of an arbitrary group $G$ with nilpotence classes $c(H)$ and $c(K)$, respectively, is a nilpotent subgroup $HK$ of nilpotence class at most $c(H)+c(K)$. We recall its consequence for pro-Lie groups. By convention, $c(1)=0$ and $c(G)=\infty$ when $G$ is not nilpotent.

\begin{theorem}[Third Main Theorem: Fitting Lemma for Pro-Lie groups]\label{t:fitting}
Let $G$ be a pro-Lie group, and let $H$ and $K$ be closed normal nilpotent subgroups of $G$. Then $L=\overline{\langle H,K\rangle}$ is a closed normal nilpotent subgroup of $G$ with $c(L)\leq c(H)+c(K)$.
\end{theorem}

\begin{proof}
Put the two nonnegative integers $r=c(H)$ and $s=c(K)$ as the nilpotence classes of $H$ and $K$, respectively. If either subgroup is trivial, the assertion is immediate. Since $H$ and $K$ are normal in $G$, their product $HK$ is a normal subgroup of $G$, and $\langle H,K\rangle=HK$. Consequently, $L=\overline{HK}$ is a closed normal subgroup of $G$. By the Fitting theorem for abstract groups \cite[Theorem 5.2.8]{robinson}, the product of two normal nilpotent subgroups of classes $r$ and $s$ is nilpotent of class at most $r+s$. Hence $\gamma_{r+s+1}(HK)=1$. Consider the commutator map
\begin{equation}\label{eq:kappa}
\begin{aligned}
\kappa_{r+s+1}\colon L^{r+s+1}&\longrightarrow L,\\
(x_1,\ldots,x_{r+s+1})&\longmapsto[x_1,\ldots,x_{r+s+1}],
\end{aligned}
\end{equation}
where the commutators are left-normed. This map is continuous, since multiplication and inversion are continuous. Since $HK$ is dense in $L$, the subset $(HK)^{r+s+1}$ is dense in $L^{r+s+1}$. The restriction of $\kappa_{r+s+1}$ to $(HK)^{r+s+1}$ is constantly equal to the identity. Therefore, by continuity and the Hausdorff property, $\kappa_{r+s+1}(L^{r+s+1})=\{1\}$. Thus $\gamma_{r+s+1}(L)=1$, and consequently $c(L)\leq r+s=c(H)+c(K)$.
\end{proof}

The argument of Theorem~\ref{t:fitting} works for every Hausdorff topological group; the pro-Lie hypothesis is needed only if one wants the resulting closed subgroup to be a pro-Lie group. It provides a bound on the nilpotence class of the closed subgroup generated by two normal nilpotent subgroups. Theorem~\ref{t:tychonovnonabelian}, on the other hand, also holds in $\ProLieNil_c$ with $c\geq3$, as observed above. For ordinary coproducts, and for conditional coproducts without the commuting hypothesis, the situation is different: distinct images in the target need not commute, so they cannot in general be obtained from a cartesian product. For instance, two identity morphisms from $\Heis(\R)$ to itself form a convergent family, but do not extend to a homomorphism $\Heis(\R)\times\Heis(\R)\to\Heis(\R)$ along the canonical embeddings. Such questions are distinct from the universal property studied here. For the existence of ordinary coproducts in categories of locally compact groups, see \cite{chir1,gm,hm77}.




\end{document}